\documentclass{amsart}

\title
{Dimers for degenerating families of M-curves}

\author{Takashi Ichikawa} 
\address{Department of Mathematics, Faculty of Science and Engineering, Saga University, 
Saga 840-8502, Japan}
\email{ichikawn@cc.saga-u.ac.jp}

\usepackage{url}
\usepackage{color}
\usepackage{rotating}

\countdef\x=23
\countdef\y=24
\countdef\z=25
\countdef\t=26

\def\tbox(#1,#2)#3{
\x=#1 \y=#2 
\multiply\x by 12 
\multiply\y by 12 
\z=\x \t=\y
\advance\z by 12 
\advance\t by 12 
\psline(\x,\y)(\x,\t)(\z,\t)(\z,\y)(\x,\y)
\advance\x by 6
\advance\y by 6 
\rput(\x,\y){{\bf #3}}}

\usepackage{amssymb}
\usepackage{graphicx}
\usepackage{epstopdf}
\usepackage{amscd}
\usepackage{appendix}
\usepackage{graphics}
\usepackage{tikz}
\def\proof{\par{\it Proof}. \ignorespaces}
\def\endproof{{\ \vbox{\hrule\hbox{%
     \vrule height1.3ex\hskip0.8ex\vrule}\hrule }}\par}

\theoremstyle{definition}

\theoremstyle{remark}

\numberwithin{equation}{section}

\let\trueint=\int
\let\truesum=\sum
\def\int{\mathop{\textstyle\trueint}\limits}
\def\sum{\mathop{\textstyle\truesum}\limits}
\let\<=\langle
\let\>=\rangle

\def\t{\mathbf{t}}
\def\0{\mathbf{0}}
\def\z{\mathbf{z}}

\def\edge{\ar@{-}}
\def\dedge{\ar@{.}}

\usepackage{amsmath, amsthm, amssymb, amsbsy,amsfonts,amscd}
\usepackage[mathscr]{eucal}
\usepackage{epsfig}
\usepackage{young}

\usepackage{amsfonts}
\usepackage{xy}
\usepackage{amssymb}
\usepackage{amsmath}

\newcommand{\thmrefer}[1]{\renewcommand\thetheorem
  {\protect\ref{#1}}\addtocounter{theorem}{-1}}

\xyoption{all}
\CompileMatrices

\begin{document}

 \begin{abstract}
We study dimer models on infinite minimal graphs with Fock's weights 
for degenerating families of M-curves of any genus 
based on works of Boutillier-Cimasoni-de Tili\`{e}re 
and Bobenko A.I.-Bobenko N.-Suris for a fixed M-curve. 
We show that these dimer models for families of M-curves behave consistently 
under their degeneration shrinking real circles to points, 
and that they can be calculated as power series in the associated deformation parameters 
which are regarded as the perturbation of Kenyon's critical dimer models. 
\end{abstract}

\maketitle

\vspace{2ex}

\vspace{2ex}


\section{Introduction}
In \cite{BouCdT2, BouCdT3}, 
Boutillier, Cimasoni and de Tili\`{e}re developed a general approach to dimer models 
on infinite minimal graphs with Fock's weights \cite{Fo}. 
Starting from M-curves which are defined as real smooth (algebraic) curves 
with maximal number (= genus\,$+\,1$) of components, 
they constructed integrable dimer models whose weights and kernels are 
are explicitly described by the associated Abelian functions, e.g., theta functions and prime forms. 
In the periodic case, 
these dimer models corresponds to Harnack curves studied by 
Kenyon, Okounkov and Sheffield \cite{KO1, KO2, KOS}. 
There were works \cite{BoB, BoBS} of Bobenko (A.I.), Bobenko (N.) and Suris 
on related results containing explicit formulas of Fock's weights and kernels 
by the Schottky uniformization theory. 

The aim of this paper is to extend the theoretical and computational results 
of \cite{BoB, BoBS, BouCdT2, BouCdT3} to degenerating M-curves of any genus 
by using variational formulas of Abelian differentials and functions. 
The main results of this paper on dimer models on infinite minimal graphs 
with Fock's weights for families of M-curves are summarized as follows: 
\begin{enumerate}

\item 
The dimer models for families of M-curves behave consistently 
under their degeneration shrinking real circles to points 
which were considered by Kenyon, Okounkov, Sheffield and Olarte \cite{KO1, KOS, O} 
for Harnack curves. 
More precisely, these models have natural limits as the dimer models 
for the normalizations of the singular curves obtained by the degeneration. 

\item 
The dimer models for families of M-curves degenerating to rational singular curves 
can be calculated as power series in the associated deformation parameters 
which are regarded as the perturbation of Kenyon's critical dimer models \cite{K}. 

\end{enumerate}

We expect that similar results hold for main objects on such dimer models 
containing amoeba maps, Ronkin functions, surface tensions and height functions 
since they are essentially expressed by Abelian functions. 
Note that the behavior of these objects under degeneration of M-curves was discussed in  \cite[Section 8.1]{BouCdT2} for the genus $1$ case and in \cite[Corollary 3]{MR} for Harnack curves, 
and seems to be important in extending results of Berggren and Borodin \cite{BerB1} 
on Aztec dimer models when the associated Harnack curves are singular 
(cf. \cite[Section 4]{BerB1} and \cite[Section 6]{BoudT}). 
See also \cite{BerB2, L} for the relationship with the tropical limits of 
imaginary normalized Abelian differentials and the crystallization of the Aztec diamonds.

\section{Abelian functions on M-curves}

\subsection{Abelian functions on an M-curve} 

An M-curve is a compact Riemann surface $R$ with an anti-holomorphic involution 
$\sigma$ whose set of fixed points is given by $g + 1$ topological circles, 
where $g$ is the genus of  $R$. 
Take a real circle $A_{0}$ in an M-curve $R$, 
and denote by $A_{1},..., A_{g}$ the remaining circles. 
Then there exist homology cycles $B_{1},.., B_{g}$ in $H_{1}(R, {\mathbb Z})$ 
such that $\sigma_{*}(B_{i}) = - B_{i}$ and that $\{ A_{1},..., A_{g}, B_{1},..., B_{g} \}$ 
gives a symplectic basis of $H_{1}(R, {\mathbb Z})$. 
Denote by $(\omega_{R})_{1},.., (\omega_{R})_{g}$ 
the Abelian differentials of the 1-st kind on $R$ which are normalized for $A_{i}$ 
in the sense that $\displaystyle \int_{A_{i}} (\omega_{R})_{j}$ is Kronecker's delta $\delta_{ij}$. 
Then as is shown in \cite{BBEIM, DuN, F}, 
the period matrix for $(R, \{ A_{i}, B_{i} \})$ defined as 
$$
\Omega_{R} = \Omega_{(R, \{ A_{i}, B_{i} \})} = 
\left( \int_{B_{i}} (\omega_{R})_{j} \right)_{1 \leq i,j \leq g} 
$$
is purely imaginary, 
and the Riemann theta function 
$$
\Theta_{R}(z) = \Theta_{(R, \{ A_{i}, B_{i} \})}(z) = \sum_{u \in {\mathbb Z}^{g}} 
\exp \left( \pi \sqrt{-1} \left( u \Omega_{R} u^{T} + 2 z u^{T} \right) \right)
$$
for $(R, \{ A_{i}, B_{i} \})$ takes real and positive values at $z \in {\mathbb R}^{g}$. 

We recall the classical definition of the prime form following \cite{Bog, F}. 
Let $(\omega_{R})_{x,y}$ denote the unique Abelian differential on $R$ 
of the $3$-rd kind whose poles are simple at $x, y$ with residues $1, -1$ respectively 
such that $\displaystyle \int_{A_{i}} (\omega_{R})_{x,y} = 0$ for all $i = 1,..., g$. 
For points $\alpha, \beta$ on an M-curve $R$ with local coordinates, 
take their lifts $\tilde{\alpha}, \tilde{\beta}$ on the universal cover $\tilde{R}$ of $R$ 
which give rise to a homotopy class of paths from $\beta$ to $\alpha$ 
and the associated line integral. 
Then the associated prime form is defined as 
$$
E_{R} (\tilde{\alpha}, \tilde{\beta}) = E_{(R, \{ A_{i} \})} (\tilde{\alpha}, \tilde{\beta}) = 
\frac{P_{R}(\tilde{\alpha}, \tilde{\beta})}{\sqrt{d \alpha} \sqrt{d \beta}}, 
$$
where 
$$
P_{R}(\tilde{\alpha}, \tilde{\beta}) = P_{(R, \{ A_{i} \})}(\tilde{\alpha}, \tilde{\beta}) := 
\left( \lim_{x \rightarrow \alpha, y \rightarrow \beta} 
\frac{-(\alpha - x)(\beta - y)}{\exp \left( \int_{\tilde{\beta}}^{\tilde{\alpha}} 
(\omega_{R})_{x, y} \right)} \right)^{1/2} 
$$
is continuous in $\alpha, \beta$ and satisfies 
$$
P_{R}(\tilde{\alpha}, \tilde{\beta}) = (\alpha - \beta)(1 + \mbox{higher order terms}) 
$$ 
for the local coordinate at $\alpha$ if $\beta$ is close to $\alpha$. 
This implies 
$$
P_{R}(\tilde{\alpha}, \tilde{\beta}) = (\alpha - \beta) \exp 
\left( - \frac{1}{2} \int_{\tilde{\beta}}^{\tilde{\alpha}} (\omega_{R})^{*}_{x, y} \right); 
\quad (\omega_{R})^{*}_{x, y} := 
(\omega_{R})_{x, y} - \left( \frac{dz}{z - x} - \frac{dz}{z - y} \right) 
\eqno(2.1)
$$
if $\alpha, \beta$ are close to each other.

\subsection{Abelian differentials on degenerating M-curves} 

We study the variation of Abelian differentials on degenerating M-curves. 
Denote by $\omega_{C/S}$ the dualizing sheaf on a nodal curve $C$ over a base scheme $S$ 
which is called the canonical invertible sheaf in \cite[Section 1]{DM}. 
The sheaf $\omega_{C/S}$ is invertible on $C$, and is functorial for $S$. 
Furthermore, 
if $S = {\rm Spec}(k)$ for an algebraically closed field $k$ and 
$\nu : C' \rightarrow C$ is the normalization of $C$ with points $p_{i}, q_{i} \in C'$ 
such that $\nu(p_{i}) = \nu(q_{i})$ are the double points on $C$, 
then $\omega_{C/S}$ is the sheaf of $1$-forms $\eta$ on $C'$ which are regular 
except for simple poles at $p_{i}, q_{i}$ such that 
${\rm Res}_{p_{i}}(\eta) + {\rm Res}_{q_{i}}(\eta) = 0$ (cf. \cite[Section 1]{DM}). 

Let $I$ be a subset of $\{ 1,..., g \}$, 
and put $J = \{ 1,..., g \} - I$. 
For a sufficiently small positive number $r$, 
we consider a family ${\mathcal R}_{I} = \{ R_{s} \}$ of stable complex curves over 
$$
U_{r} = \left\{ \left. s = (s_{i})_{i \in I} \ \right| \, 0 \leq s_{i} < r \right\} 
$$ 
such that for $s_{i} > 0$ $(i \in I)$, 
$R_{s}$ are M-curves with real circles $A_{0}, A_{1},..., A_{g}$ as above and that $s_{i} \downarrow 0$  corresponds to the shrinking of $A_{i}$ to a point $\alpha_{i}$ (See Figure 1). 
\mbox{}

\vspace{-5ex}
\hspace{5ex}  \includegraphics[width=150mm]{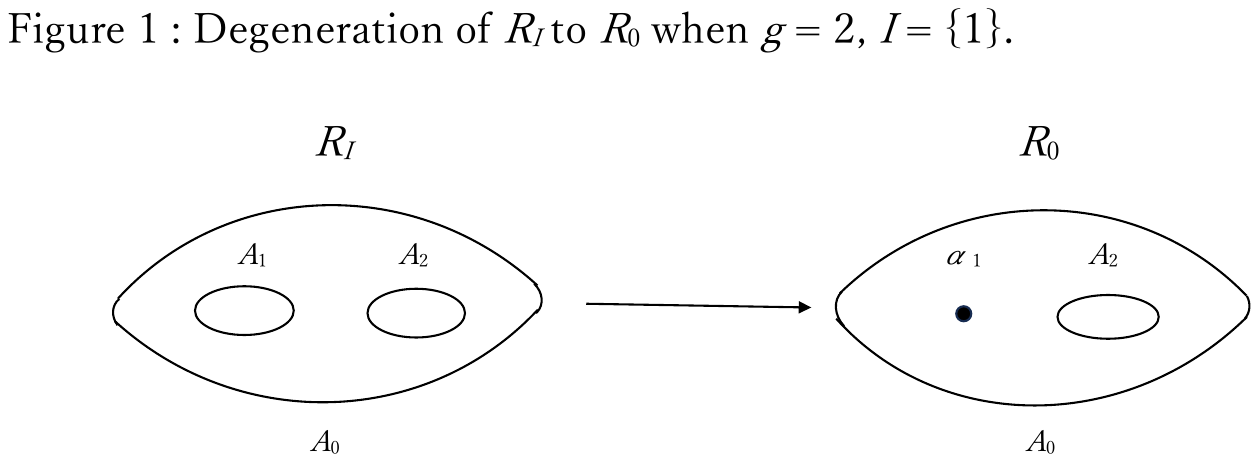}

\vspace{-30ex}
\noindent
Then the normalization $R_{J}$ of the stable curve $R_{0} = R_{s}|_{s = 0}$ 
becomes an M-curve with real circles $A_{0}$, $A_{j}$ $(j \in J)$ and 
the anti-holomorphic involution $\sigma$, 
and $R_{0}$ is obtained from $R_{J}$ by identifying points 
$\alpha_{i}$ and $\sigma(\alpha_{i})$ $(i \in I)$ on $R_{J}$ in pairs. 
Therefore, $\{ R_{s} \}$ gives the deformation of $R_{0}$ by the relation 
$$
\xi_{\alpha_{i}} \cdot \xi_{\sigma(\alpha_{i})} = s_{i} \quad (i \in I)
$$
for certain local coordinates $\xi_{p}$ at $p \in R_{0}$ satisfying $\xi_{p}(p) = 0$. 
Such families can be constructed by Schottky uniformization theory of M-curves 
\cite[Section 10]{BoBS} (see Section 2.3 below). 
We call global sections of the dualizing sheaf $\omega_{{\mathcal R}/U_{r}}$ 
{\it stable Abelian differentials} on ${\mathcal R}/U_{r}$. 
\vspace{2ex}

\noindent 
{\bf Theorem 2.1.} 
\begin{it}
There exists a unique basis $\{ \omega_{i} \}_{1 \leq i \leq g}$ 
of stable Abelian differentials on ${\mathcal R}_{I}/U_{r}$ 
which is normalized for $\left\{ A_{i} \right\}_{1 \leq i \leq g}$ when all $s_{i}$ are positive. 
Furthermore, 
the pullback $(\omega_{j}|_{R_{0}})^{*}$ of $\omega_{j}|_{R_{0}}$ to $R_{J}$ $(j \in J)$ are 
the Abelian differentials on $R_{J}$ of the 1st kind normalized for $A_{j}$. 
\end{it}
\vspace{2ex}

\noindent 
{\it Proof.} 
The assertion was analytically proved in \cite{F} for one step degenerations, 
and was algebro-geometrically proved in \cite[Theorem 7.1]{I3} for general degenerations. 
\ $\square$ 
\vspace{2ex}

\noindent
{\bf Theorem 2.2.} 
\begin{it} 
Let $x, y$ be non-intersecting smooth sections from $U_{r}$ into smooth points on ${\mathcal R}_{I}$. 
Then there exists a unique meromorphic section $\omega_{x, y}$ of $\omega_{{\mathcal R}_{I}/U_{r}}$ 
such that  
\begin{itemize}

\item 
$\omega_{x, y}$ is regular outside $x, y$, 
and it has simple poles at $x, y$ with residues $1, -1$ respectively, 

\item 
When all $s_{i}$ $(i \in I)$ are positive, 
$\displaystyle \int_{A_{i}} \omega_{x, y} = 0$ on $R_{s}$ for all $A_{i}$ $(i =1,..., g)$. 

\end{itemize}
Furthermore, 
the pullback $(\omega_{x, y}|_{R_{0}})^{*}$ of $\omega_{x, y}|_{R_{0}}$ to $R_{J}$ satisfies 
$\displaystyle \int_{A_{j}} (\omega_{x, y}|_{R_{0}})^{*} = 0$ for all $A_{j}$ $(j \in J)$. 
\end{it}
\vspace{2ex}

\noindent
{\it Proof.} 
By the above description of $\omega_{C/S}$,  
there exists a unique meromorphic Abelian differential $(\omega_{x, y})_{0}$ on $R_{0}$ 
such that its pullback $((\omega_{x, y})_{0})^{*}$ to $R_{J}$ satisfies 
$\displaystyle \int_{A_{j}} ((\omega_{x, y})_{0})^{*} = 0$ for all $A_{j}$ $(j \in J)$ 
and has only poles at $x, y$ which are simple with residue $1, -1$ respectively. 
Then the Riemann-Roch theorem and residue theorem for nodal curves stated in 
\cite[Section A of Chapter 3]{HM} implies that 
the coherent sheaf of meromorphic sections of $\omega_{{\mathcal R}_{I}/U_{r}}$ 
whose poles are all simple located at $x, y$ is locally free on $U_{r}$ of rank $g + 1$. 
Therefore, 
$(\omega_{x, y})_{0}$ can be extended to a meromorphic section 
$(\omega_{x, y})_{s}$ of $\omega_{{\mathcal R}_{I}/U_{r}}$ with simple poles 
at $x, y$ with residues $1, -1$ respectively by retaking $r > 0$ if necessary. 
Then by Theorem 2.1, 
$$
(\omega_{x, y})_{s} - \sum_{i=1}^{g} \int_{A_{i}} (\omega_{x, y})_{s} \cdot \omega_{i}
$$
gives the above $\omega_{x, y}$. 
\ $\square$

\subsection{Variation of Abelian functions} 

Let the notation be as in Section 2.2. 
Then for each $s = (s_{i})_{i \in I} \in U_{r}$ with $s_{i} > 0$ $(i \in I)$, 
each period matrix 
$$
\Omega_{R_{s}} = \Omega_{\left( R_{s}, \{ A_{i}, B_{i} \}_{1 \leq i \leq g} \right)} 
$$
gives a symmetric bilinear form on ${\mathbb Z}^{g} \times {\mathbb Z}^{g}$, 
where 
$$
{\mathbb Z}^{g} = \bigoplus_{i = 1}^{g} {\mathbb Z} B_{i} 
= {\mathbb Z}^{I} \oplus {\mathbb Z}^{J}; \quad 
{\mathbb Z}^{I} := \bigoplus_{i \in I} {\mathbb Z} B_{i}, \ 
{\mathbb Z}^{J} := \bigoplus_{j \in J} {\mathbb Z} B_{j}, 
\eqno(2.2)
$$
and $\{ A_{j}, B_{j} \}_{j \in J}$ gives a symplectic basis of $R_{J}$. 

In the following theorem, 
$s \downarrow 0$ means that $s_{i} \downarrow 0$ for all $i \in I$. 
\vspace{2ex}

\noindent
{\bf Theorem 2.3.}  
\begin{it} 

{\rm (1)} 
For $z = (z_{I}, z_{J}) \in {\mathbb C}^{g} = {\mathbb C}^{I} \oplus {\mathbb C}^{J}$, 
$$
\lim_{s \downarrow 0} \Theta_{(R_{s}, \{ A_{i}, B_{i} \}_{1 \leq i \leq g})} (z) = 
 \Theta_{(R_{J}, \{ A_{j}, B_{j} \}_{j \in J})} (z_{J}). 
$$ 

{\rm (2)} 
For non-intersecting sections  $\alpha, \beta$ from $U_{r}$ into smooth points 
on ${\mathcal R}_{I}$,
and their smooth lifts $\tilde{\alpha}, \tilde{\beta}$ respectively 
to the universal covers of $R_{s}$ with $s_{i} > 0$ $(i \in I)$, 
$$
\lim_{s \downarrow 0} E_{(R_{s}, \{ A_{i} \}_{1 \leq i \leq g})} (\tilde{\alpha}, \tilde{\beta}) 
= E_{(R_{J}, \{ A_{j} \}_{j \in J})} (\tilde{\alpha}, \tilde{\beta}), 
$$ 
where $\tilde{\alpha}, \tilde{\beta}$ in the right hand side denotes the projection to 
the universal cover of $R_{J}$ associated with the surjective group homomorphism 
$\pi_{1}(R_{s}) \rightarrow \pi_{1}(R_{J})$ sending $A_{i}, B_{i}$ $(i \in I)$ to the identity 
and  $A_{j}, B_{j}$ $(j \in J)$ to $A_{j}, B_{j}$ respectively. 
\end{it} 
\vspace{2ex}

\noindent 
{\it Proof.} 
The assertion (1) follows from the variational formula shown in \cite{F, Go, M} 
that for each $i \in I$, 
the $(i, i)$-component $(\Omega_{R_{s}})_{ii}$ of $\Omega_{R_{s}}$ satisfies 
$$
\lim_{s_{i} \downarrow 0} \exp \left( 2 \pi \sqrt{-1} (\Omega_{R_{s}})_{ii} \right) = 0. 
$$
The assertion (2) follows from (2.1) and Theorem 2.2. 
\ $\square$

\subsection{Schottky uniformization} 

Assume that $I = \{ 1,..., g \}$ and $J = \emptyset$. 
Then as is shown in \cite[Section 10]{BoBS}, 
we can construct familes of M-curves ${\mathcal R}_{I} = \{ R_{s} \}$ of genus $g$ 
with $\sigma(z) = \overline{z}$ which are uniformized by the Schottky groups $\Gamma$ 
generated by $\gamma_{i} \in PGL_{2}({\mathbb C})$ $(i \in I)$, 
where 
$$
\frac{\gamma_{i}(z) - \alpha_{i}}{\gamma_{i}(z) - \overline{\alpha_{i}}} = 
s_{i} \frac{z - \alpha_{i}}{z - \overline{\alpha_{i}}} 
\quad \left( z \in {\mathbb P}^{1}_{\mathbb C} \right)   
$$
for complex numbers $\alpha_{i}$ with positive imaginary part and $0 < s_{1},..., s_{g} < 1$ 
(See Figure 2). 

\vspace{-5ex}
\hspace{5ex}  \includegraphics[width=150mm]{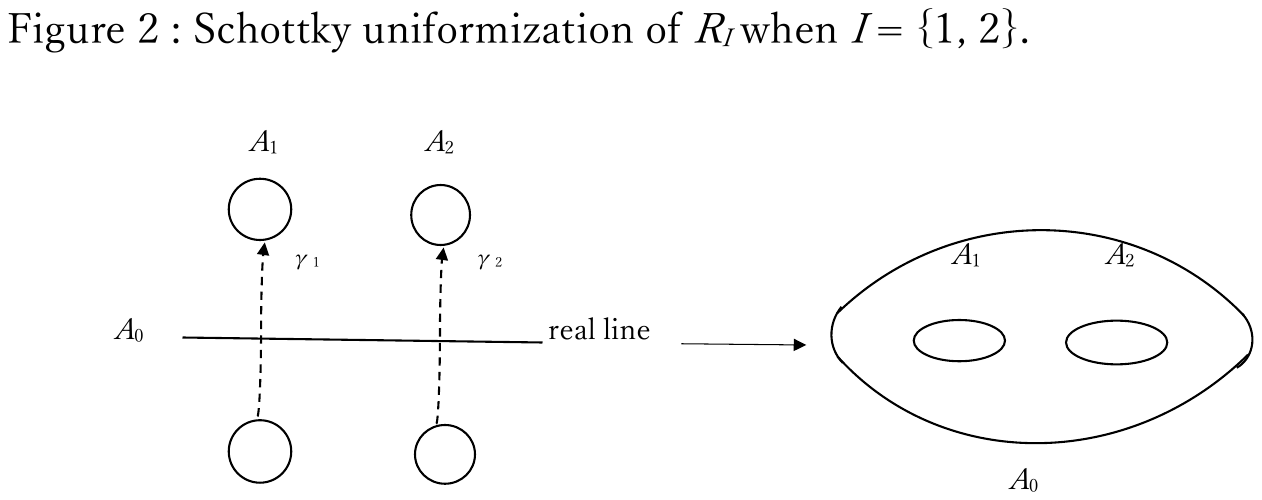}

\vspace{-27ex}
\noindent
Therefore, 
${\mathcal R}_{I}$ gives the deformation by $s_{i}$ of the singular curve $R_{0}$ obtained by 
identifying $\alpha_{i} = \overline{\alpha_{i}} \in {\mathbb P}^{1}_{\mathbb C}$ $(i \in I)$ 
for the local coordinates 
$$
\xi_{\alpha_{i}}(z) = \frac{z - \alpha_{i}}{z - \overline{\alpha_{i}}}, 
\quad \xi_{\overline{\alpha_{i}}}(z) = \frac{z - \overline{\alpha_{i}}}{z - \alpha_{i}}. 
$$
By results of \cite{BoB, BoBS} and \cite{I1, I2}, 
the above Abelian differentials $(\omega_{R_{s}})_{i}$, $(\omega_{R_{s}})_{x, y}$ 
on $R_{s}$ of the 1st, 3rd kind respectively, 
and the associated multiplicative periods $\exp(2 \pi \sqrt{-1} (\Omega_{R_{s}})_{ij})$ 
are expressed as 
\begin{eqnarray*}
(\omega_{R_{s}})_{i} & = & 
\frac{1}{2 \pi \sqrt{-1}} 
\sum_{\gamma \in \Gamma/\langle \gamma_{i} \rangle} 
\left( \frac{1}{z - \gamma(\alpha_{i})} - \frac{1}{z - \gamma(\overline{\alpha_{i}})} \right) dz, 
\\
(\omega_{R_{s}})_{x, y} & = & 
\sum_{\gamma \in \Gamma} 
\left( \frac{1}{z - \gamma(x)} - \frac{1}{z - \gamma(y)} \right) d z, 
\\ 
\exp(2 \pi \sqrt{-1} (\Omega_{R_{s}})_{ij}) & = & 
\prod_{\gamma \in \langle \gamma_{i} \rangle \backslash \Gamma / \langle \gamma_{j} \rangle} 
\psi_{ij} (\gamma), 
\end{eqnarray*} 
where 
$$
\psi_{ij} = \left\{ \begin{array}{ll} 
s_{i} 
& (\mbox{$i = j$ and $\gamma \in \langle \gamma_{i} \rangle$}), 
\\ 
\left[ \alpha_{i}, \overline{\alpha_{i}}; \gamma(\alpha_{j}), \gamma(\overline{\alpha_{j}}) \right] 
& (\mbox{otherwise}) 
\end{array} \right. 
$$
with $[a, b; c, d]$ denotes the cross-ratio: 
$$
\frac{(a-c)(b-d)}{(a-d)(b-c)}. 
$$
In what follows, 
we assume that these infinite sums and products are absolutely convergent 
which condition is discussed in \cite{BBEIM}, 
especially, is guaranteed when $|s_{1}|,..., |s_{g}|$ are sufficiently small. 
Furthermore, by \cite[Proposition 3.2]{I1} and its proof, 
especially the formula 
$$
\gamma_{\pm i}(v) - \gamma_{\pm i}(w) = 
\frac{(\alpha_{i} - \alpha_{-i})^{2} (v - w) s_{i}}
{(v - \alpha_{\mp i} - s_{i}(v - \alpha_{\pm i}))(w - \alpha_{\mp i} - s_{i}(w - \alpha_{\pm i}))}; 
\ \ \gamma_{-i} := \gamma_{i}^{-1}, \ \alpha_{-i} := \overline{\alpha_{i}}, 
$$  
we can express these objects by computable power series in $s_{1},..., s_{g}$ as follows: 
\begin{eqnarray*}
(\omega_{R_{s}})_{i} & = & 
\frac{1}{2 \pi \sqrt{-1}} 
\left( \frac{1}{z - \alpha_{i}} - \frac{1}{z - \overline{\alpha_{i}}} 
+ O(\text{max}_{k} \{ s_{k} \}) \right) dz, 
\\
(\omega_{R_{s}})_{x, y} & = & 
\left( \frac{1}{z - x} - \frac{1}{z - y} + O(\text{max}_{k} \{ s_{k} \}) \right) dz, 
\\ 
\exp(2 \pi \sqrt{-1} (\Omega_{R_{s}})_{ij}) & = & 
\left\{ \begin{array}{ll} 
s_{i} + O(\text{max}_{k,l} \{ s_{k} s_{l} \}) 
& (\mbox{$i = j$ and $\gamma \in \langle \gamma_{i} \rangle$}), 
\\ 
\left[ \alpha_{i}, \overline{\alpha_{i}}; \alpha_{j}, \overline{\alpha_{j}} \right] + 
O(\text{max}_{k} \{ s_{k} \}) 
& (\mbox{otherwise}). 
\end{array} \right. 
\end{eqnarray*} 
Therefore, 
$\Theta_{(R_{s}, \{ A_{i}, B_{i} \})}(z)$ can be expressed as a computable power series 
$$
\Theta_{(R_{s}, \{ A_{i}, B_{i} \})}(z) = 1 + 
\sum_{i=1}^{g} \left( e^{2 \pi \sqrt{-1} z_{i}} + e^{-2 \pi \sqrt{-1} z_{i}} \right) \sqrt{s_{i}} 
+ O(\text{max}_{k,l} \{ \sqrt{s_{k}} \sqrt{s_{l}} \}) 
$$
in $\sqrt{s_{1}},..., \sqrt{s_{g}}$ with constant term $1$. 
Furthermore, as is shown in \cite{Bog}, 
$P_{R_{s}}(\ast, \ast)$ is defined for $\alpha, \beta$ in the domain of discontinuity of $\Gamma$, 
and 
$$
P_{R_{s}}(\alpha, \beta) 
= (\alpha - \beta) \prod^{*}_{\gamma \in \Gamma - \{ 1 \}} 
\left[ \alpha, \beta; \gamma(\beta), \gamma(\alpha) \right] 
= (\alpha - \beta) \prod^{*}_{\gamma \in \Gamma - \{ 1 \}} 
\left( 1 - \left[ \alpha, \gamma(\beta); \beta, \gamma(\alpha) \right] \right), 
$$
where the product $\displaystyle \prod^{*}_{\gamma \in \Gamma - \{ 1 \}}$ 
means that only one element in each pair $\gamma, \gamma^{-1}$ should be taken. 
Therefore, 
$P_{R_{s}}(\alpha, \beta)$ is also expressed as power series in $s_{1},..., s_{g}$ as 
$$
P_{R_{s}}(\alpha, \beta) = (\alpha - \beta) 
\left( 1 + \sum_{i=1}^{g} 
[\alpha, \alpha_{i}; \beta, \overline{\alpha_{i}}] \cdot 
[\alpha, \overline{\alpha_{i}}; \beta, \alpha_{i}] s_{i} \right) 
+ O(\text{max}_{k,l} \{ s_{k} s_{l} \}) 
$$
with constant term $\alpha - \beta$.

\section{Dimers for degenerating M-curves}

\subsection{Dimers for an M-curve}

Following \cite[Sections 3.1 and 3.2]{BouCdT3}, 
we review dimer models for smooth M-curves. 
Let $G = (V, E)$ be an infinite and locally finite graph embedded in the plane with faces 
given as topological discs. 
If $G^{*} = (V^{*}, E^{*})$ denotes the dual embedded graph, 
then the associated quad-graph $G^{\diamond}$ is obtained from the vertex set 
$V \sqcup V^{*}$ by joining a vertex $v \in V$ and a dual vertex $f \in V^{*}$ 
with an edge when $v$ lies on the boundary of the face corresponding to $f$. 
Then $G^{\diamond}$ embeds in the plane with (possibly degenerate) 
quadrilaterals faces (see Figure 3). 
We define a train-track of $G$ as a maximal chain of adjacent quadrilaterals 
of $G^{\diamond}$ such that when one enters a quadrilateral, 
one exits through the opposite edge. 
We assume that $G$ is bipartite, 
namely that $V$ has a partition $B \sqcup W$ into black and white vertices 
such that no edge of $E$ connects two vertices of the same color, 
and then train-tracks can be consistently oriented with black vertices on the right
and white vertices on the left of the path (see Figure 3). 
Denote by ${\mathcal T}$ the set of consistently oriented train-tracks of 
the bipartite graph $G$, 
and a bipartite and planar graph $G$ is said to be minimal \cite{GK} 
if its train-tracks do not self-intersect, 
and no pair of oriented train-tracks intersect twice in the same direction. 
Then train tracks do not form loops, 
and that $G$ has neither multiple edges, nor degree $1$ vertices, 
and hence a minimal graph is a simple graph. 

\vspace{-5ex}
\hspace{5ex}  \includegraphics[width=150mm]{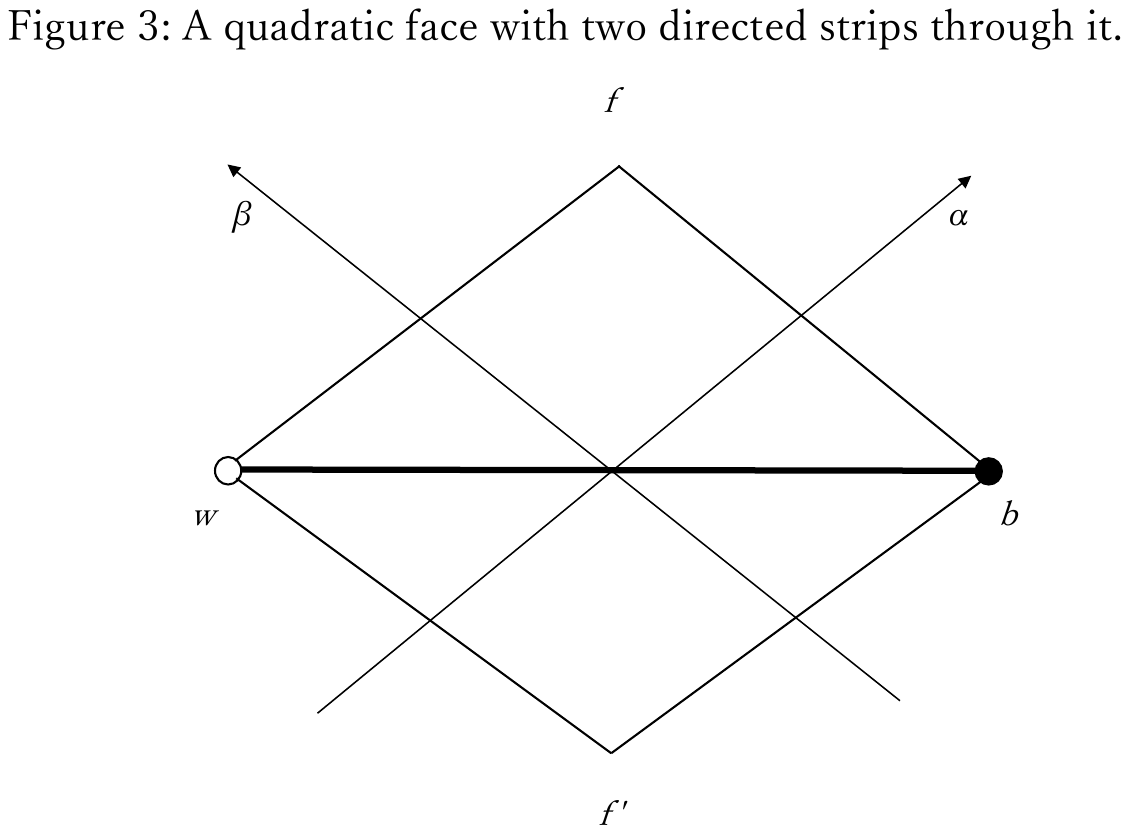}

\vspace{-10ex}
Let $R$ be an M-curve. 
Then it has an oriented real component denoted by $A_{0}$ 
containing the base point $x_{0}$. 
We associate each oriented train-track $T \in {\mathcal T}$ of $G$ 
with an element $\alpha_{T}$ of $A_{0}$ called its angle 
since $A_{0}$ is naturally identified with ${\mathbb R}/{\mathbb Z}$. 
Following Fock \cite{Fo}, 
we define a function $d$ from the set of vertices of $G^{\diamond}$ 
into the divisor group ${\rm Div}(R)$ of $R$ as follows. 
Fix a face $f_{0}$ and put $d(f_{0}) = 0$. 
Then along an edge of $G^{\diamond}$ crossing a train-track $T$ with angle $\alpha_{T}$, 
we add $\alpha_{T}$ to the value of $d$ if we arrive at a black vertex or leave a white vertex, 
and hence in Figure 3, 
$$
d(b) = d(f) + \alpha, \quad d(w) = d(f) - \beta, \quad d(f') = d(f) + \alpha - \beta
$$
for the associated angles $\alpha, \beta \in A_{0}$. 
Therefore, the degree of $d(x) \in {\rm Div}(R)$ is equal to 
$$
{\rm deg} \, d(x) = \left\{ \begin{array}{ll} 
1 & \mbox{if $x$ is a black vertex of $G$,} 
\\ 
0 & \mbox{if $x$ is a face of $G$,} 
\\ 
-1 & \mbox{if $x$ is a white vertex of $G$,} 
\end{array} \right. 
$$
and hence, for any face $f$ of $G$, 
$d(f)$ belongs to the subgroup ${\rm Div}^{0}(R) \subset {\rm Div}(R)$ of divisors on $R$ 
with degree $0$. 
Then its image by the Abel-Jacobi map: 
$$
d(f) \mapsto \sum_{k} \int_{\alpha_{k}}^{\beta_{k}} ((\omega_{R})_{1},..., (\omega_{R})_{g}) 
\ \ \text{for} \ \ d(f) = \sum_{k} (\beta_{k} - \alpha_{k}) 
$$
belongs to $({\mathbb R}/{\mathbb Z})^{g}$ 
which is contained in the Jacobian variety $J_{R} = {\rm Jac}(R)$ of $R$. 
Note that only a special class of angle assignments $T \mapsto \alpha_{T}$ 
gives rise to probabilistic models. 
It can be described as follows, see \cite{BouCdT1} for more detail. 
We call two non-closed oriented planar curves parallel (resp. antiparallel) 
when they intersect infinitely many times in the same direction (resp. in opposite directions), 
or when they are disjoint and cross a topological disc in the same direction 
(resp. in opposite directions). 
We consider a triple of oriented train-tracks of $G$, pairwise non-parallel. 
Let $D$ be a compact disk outside of which these train-tracks do not meet, 
apart from possible anti-parallel ones, 
and order this triple of elements of ${\mathcal T}$ cyclically according to 
the outgoing points of the corresponding oriented curves in the circle $\partial D$. 
Then we have a well-defined partial cyclic order on ${\mathcal T}$ 
(cf. \cite[Section 2.3]{BouCdT1}), 
and note that $A_{0}$ is an oriented topological circle endowed with 
a total cyclic order as well which allows for the following definition. 
We define $X_{G}$ as the set of maps $\alpha: {\mathcal T} \rightarrow A_{0}$ 
that are monotone, in the sense that they preserve the cyclic order, 
and such that non-parallel train-tracks have distinct images. 
One of the main results of \cite{BouCdT1} states that if $G$ is minimal, 
then $X_{G}$ is included in the space of minimal immersions of $G$, 
and coincides with it if $G$ is minimal and periodic 
(cf. \cite[Theorem 23 and Corollary 29]{BouCdT1}).

\subsection{Fock's Kasteleyn operator with kernels and inverse}

In order to define a version of Fock's adjacency operator satisfying Kasteleyn's condition, 
fix a minimal graph $G$, an M-curve $R$ of genus $g$ with real circles $A_{i}$ $(i = 0, 1,..., g)$ 
and an angle map $\alpha \in X_{G}$. 
Fix an arbitrary lift $\tilde{\alpha} : {\mathcal T} \rightarrow \tilde{A}_{0}$ 
of $\alpha$, namely, lifts $\tilde{\alpha}_{T} \in \tilde{A}_{0}$ of the angles $\alpha_{T} \in A_{0}$, 
where $\tilde{A}_{0} \subset \tilde{R}$ denotes the universal cover of $A_{0}$. 
The Abel-Jacobi map defines an embedding of $A_{0}$ 
in $({\mathbb R}/{\mathbb Z})^{g}$, 
and hence an embedding of $\tilde{A}_{0}$ in ${\mathbb R}^{g} \subset {\mathbb C}^{g}$ 
which is identified with the universal cover $\tilde{J}_{R}$ of $J_{R}$. 
We define a lift $\tilde{d} : V(G^{\diamond}) \rightarrow {\rm Div} (\tilde{R})$ 
of the discrete Abel map $d$ by putting $\tilde{d}(f_{0}) = 0$, 
and compute the values at every vertex iteratively by adding and subtracting 
the lifts $\tilde{\alpha}$ of the crossed train-tracks with the same local rule as $d$. 
Therefore, 
if $b$ (resp. $w$) and $f$ are separated by a train-track 
with angle $\alpha$ (resp. $\beta$) as in Figure 3, 
we have $\tilde{d}(b) = \tilde{d}(f) + \tilde{\alpha}$ 
(resp. $\tilde{d}(w) = \tilde{d}(f) - \tilde{\beta}$). 
For any face $f$ of $G$, the divisor $\tilde{d}(f)$ has degree $0$, 
and its image by the Abel-Jacobi map 
${\rm Div}^{0} (\tilde{R}) \rightarrow {\mathbb C}^{g}$ belongs ${\mathbb R}^{g}$ 
which we denote by the same notation $\tilde{d}(f)$. 

Following \cite[Sections 3.3 and 3.4]{BouCdT3}, 
we review the definition of {\it Fock's adjacency operator}, its kernels and inverse. 
Fock's adjacency operator $K_{R}$ is the weighted adjacency operator of the graph $G$ 
which is indexed by elements $t \in ({\mathbb R}/{\mathbb Z})^{g} \subset J_{R}$ 
with coefficients called {\it Fock's weights} given as follows: 
for every edge $wb$ crossed by train-tracks with angles $\alpha, \beta$ in $A_{0}$ 
as in Figure 3, 
we define 
$$
(K_{R})_{w, b} = 
\frac{E_{R} (\tilde{\beta}, \tilde{\alpha})}
{\Theta_{R} (\tilde{t} + \tilde{d}(f)) \, \Theta_{R} (\tilde{t} + \tilde{d}(f'))}; \ \ 
\left\{ \begin{array}{l} 
E_{R}(*, *) = E_{(R, \{ A_{i} \})}(*, *), 
\\
\Theta_{R}(*) = \Theta_{(R, \{ A_{i}, B_{i} \})}(*), 
\end{array} \right. 
\eqno(3.1) 
$$   
where $\tilde{t} \in {\mathbb R}^{g}$ is a lift of $t \in ({\mathbb R}/{\mathbb Z})^{g}$. 
For arbitrary vertices $x, y$ of the quad-graph $G^{\diamond}$, 
a meromorphic form $(g_{R})_{x, y}$ on $R$ is defined as an element in the kernel of $K_{R}$. 
When $f$ is a vertex of $G^{*}$ and $w, b$ are white, black respectively adjacent vertices 
in $G$, 
put 
$$
\begin{array}{l}
{\displaystyle (g_{R})_{f, w}(\tilde{u}) = (g_{R})_{w, f}(\tilde{u})^{-1} = 
\frac{\Theta_{R} (\tilde{t} + (\tilde{u} + \tilde{d}(w)))}{E_{R} (\tilde{u}, \tilde{\beta})}}, 
\\ 
\\
{\displaystyle (g_{R})_{b, f}(\tilde{u}) = (g_{R})_{f, b}(\tilde{u})^{-1} = 
\frac{\Theta_{R} (-\tilde{t} + (\tilde{u} - \tilde{d}(b)))}{E_{R} (\tilde{u}, \tilde{\alpha})}} 
\end{array} 
\eqno(3.2) 
$$ 
for $\tilde{u} \in \tilde{R}$, 
where $\tilde{u} + \tilde{d}(w)$ and $\tilde{u} - \tilde{d}(b)$ are divisors on $\tilde{R}$ of degree $0$ 
and hence they can be regarded in ${\mathbb C}^{g}$ by the Abel-Jacobi map $\tilde{d}$. 
When $x, y$ are not necessarily neighbors, 
take a path $x = x_{1},...,x_{n} = y$ in $G^{\diamond}$, 
and put 
$$
(g_{R})_{x, y}(\tilde{u}) = \prod_{j=1}^{n-1} (g_{R})_{x_{j}, x_{j+1}} (\tilde{u}) 
\eqno(3.3) 
$$
which is well-defined, i.e., is independent of the choice of the path. 
Then it is shown in \cite[Proposition 31 and Lemma 33]{BouCdT3} (see also \cite{Fo}) 
that Fock's adjacency operator $K_{R}$ is a Kasteleyn operator, i.e., 
defines a dimer model on $G$. 
Furthermore, for any vertex $x$ of $G^{\diamond}$, 
$((g_{R})_{b, x}(\tilde{u}))_{b \in B}$ (resp. $((g_{R})_{x, w}(\tilde{u}))_{w \in W}$) 
are in the right (resp. left) kernels of $K$ which mean that 
$$
\sum_{b \sim w} (K_{R})_{w, b} \, (g_{R})_{b, x} (\tilde{u}) = 0 \ (w \in W), \ \ 
(\mbox{resp.}   
\sum_{w \sim b} (g_{R})_{x, w} (\tilde{u}) \, (K_{R})_{w, b} = 0 \ (b \in B)), 
$$ 
where $b \sim w, w \sim b$ denote that $b$ and $w$ are adjacent. 
Furthermore, it is shown in \cite[Theorem 40]{BouCdT3} and its proof that 
for the component $R^{+}$ of $R - \bigcup_{i=0}^{g} A_{i}$ whose oriented boundary contains $A_{0}$, 
$$
(A_{R})^{u_{0}}_{b,w} = \frac{1}{2 \pi \sqrt{-1}} \int_{C^{u_{0}}_{b,w}} (g_{R})_{b,w} 
\quad \left( u_{0} \in R^{+} - \alpha({\mathcal T}) \right) 
\eqno(3.4) 
$$
gives an inverse of the operator $K_{R}$, 
where $C^{u_{0}}_{b,w}$ is an oriented simple path from $\sigma(u_{0})$ to $u_{0}$ 
described in \cite[Section 3.5]{BouCdT3}, 
and 
$$
(K_{R})_{w,b} \cdot (g_{R})_{b,w} = (\omega_{R})_{\beta, \alpha} + \sum_{j=1}^{g} 
\left( \frac{\partial \log \Theta_{R}}{\partial z_{j}} (\tilde{t} + \tilde{d}(f)) - 
\frac{\partial \log \Theta_{R}}{\partial z_{j}} (\tilde{t} + \tilde{d}(f')) \right) (\omega_{R})_{j}. 
\eqno(3.5) 
$$  

As is explained in \cite[Section 8.1]{BouCdT2} and \cite[Section 2.3]{BoudT}, 
if the genus of $R$ is $0$, 
then $K_{w, b} = \beta - \alpha$ which is equivalent to Kenyon's critical dimer models 
\cite{K} in which the anti-holomorphic involution is defined as $\sigma(z) = 1/ \overline{z}$, 
and $A_{0} = \{ z \in {\mathbb C} \ | \ |z| = 1 \}$.

\subsection{Dimers for degenerating M-curves}

We state main results of this paper on dimer models for the family 
${\mathcal R}_{I} = \{ R_{s} \}$ of degenerating M-curves with real ovals $A_{i}$ $(i = 1,..., g)$ 
given in Section 2.2. 
Let $G$ be an infinite minimal graph, 
$\tilde{\alpha}_{s}: {\mathcal T} \rightarrow \tilde{A}_{0}$ $(s \in U_{r})$ 
be a continuous family of lifted angle maps for ${\mathcal R}_{I}$, 
and denote by $\tilde{d}: V(G^{\diamond}) \rightarrow {\rm Div}(\tilde{R}_{s})$ 
the associated lifted discrete Abel map. 
Then for each face $f$ of $G$, 
the image of $\tilde{d}(f) \in {\rm Div}^{0}(\tilde{R}_{s})$ by the Abel-Jacobi map 
${\rm Div}^{0}(\tilde{R}_{s}) \rightarrow {\mathbb R}^{g}$ is a continuous function of $s$.     
\vspace{2ex}

\noindent
{\bf Theorem 3.1.} 
\begin{it} 
Let the notation be as above. 
Then the associated Fock's weights $(K_{{\mathcal R}_{I}})_{w, b}$, 
the kernel forms $(g_{{\mathcal R}_{I}})_{x, y}(\tilde{u})$ satisfy 
$$
\lim_{s \downarrow 0} (K_{{\mathcal R}_{I}})_{w, b} = (K_{R_{J}})_{w, b}, \quad 
\lim_{s \downarrow 0} (g_{{\mathcal R}_{I}})_{x, y}(\tilde{u}) = (g_{R_{J}})_{x, y}(\tilde{u}), \quad 
$$ 
where $s \downarrow 0$ means that $s_{i} \downarrow 0$ for all $i \in I$. 
Furthermore, the inverse $(A_{{\mathcal R}_{I}})^{u_{0}}_{b,w}$ of $K_{{\mathcal R}_{I}}$ satisfies 
$$
\lim_{s \downarrow 0} (A_{{\mathcal R}_{I}})^{u_{0}}_{b, w} = (A_{R_{J}})^{u_{0}}_{b, w} 
$$
if the path $C^{u_{0}}_{b,w}$ is defined for $R_{s}$ and $R_{J}$. 
In particular, if $I = \{ 1,..., g \}$, 
then the dimer models for ${\mathcal R}_{I}$ become Kenyon's critical dimer models 
under $s \downarrow 0$. 
\end{it} 
\vspace{2ex}

\noindent
{\it Proof.} 
For each face $f$ of $G$ and $\tilde{t} \in {\mathbb R}^{g}$, 
by Theorem 2.3 (1), 
$$
\lim_{s \downarrow 0} \Theta_{(R_{s}, \{ A_{i}, B_{i} \})} (\tilde{t} + \tilde{d}(f)) = 
\Theta_{(R_{J}, \{ A_{j}, B_{j} \}_{j \in J})} (\tilde{t}_{J} + \tilde{d}(f)_{J}), 
$$ 
where $z_{J}$ is the image of $z$ by the natural projection 
${\mathbb R}^{g} = {\mathbb R}^{I} \oplus {\mathbb R}^{J} \rightarrow {\mathbb R}^{J}$ 
corresponding to (2.2). 
By Theorem 2.1, 
$\tilde{d}(f)_{J}$ is the image of $f$ by the discrete Abel map 
$\tilde{d}: V(G^{\diamond}) \rightarrow {\rm Div}(\tilde{R}_{J})$ and 
the Abel-Jacobi map ${\rm Div}^{0}(\tilde{R}_{J}) \rightarrow {\mathbb C}^{J}$ 
which are associated with the universal cover $\tilde{R}_{J}$ of $R_{J}$. 
Therefore, by Theorem 2.3 (2) and the formula (3.1), we have 
$$
\lim_{s \downarrow 0} (K_{{\mathcal R}_{I}})_{w, b} = (K_{R_{J}})_{w, b}. 
$$
Similarly, the remaining assertions can be shown by applying Theorems 2.1--2.3 to the formulas (3.3)--(3.5). 
\ $\square$ 
\vspace{2ex}

\noindent
{\bf Theorem 3.2.} 
\begin{it} 
Assume that $I = \{ 1,..., g \}$ in Theorem 3.1. 
Then the associated dimers models with Fock's weights $(K_{{\mathcal R}_{I}})_{w, b}$, 
the kernel forms $(g_{{\mathcal R}_{I}})_{x, y}(\tilde{u})$ 
and the inverse $(A_{{\mathcal R}_{I}})^{u_{0}}_{b,w}$ of $K_{{\mathcal R}_{I}}$ are expressed 
as the perturbation by the parameters $\sqrt{s_{i}}$ $(i \in I)$ 
of Kenyon's critical dimer models \cite{K}. 
\end{it}
\vspace{2ex}

\noindent
{\it Proof.} 
In the notation in Section 2.4, 
the image of $\sum_{k} (y_{k} - x_{k})$ $(x_{k}, y_{k} \in {\mathbb R})$ by the Abel-Jacobi map satisfies 
$$
\left( \exp \left( 2 \pi \sqrt{-1} \sum_{k} \int_{x_{k}}^{y_{k}} (\omega_{R_{s}})_{i} 
\right) \right)_{1 \leq i \leq g} 
= \left( \prod_{k} \prod_{\gamma \in \Gamma/\langle \gamma_{i} \rangle} 
[y_{k}, x_{k}; \gamma(\alpha_{i}), \gamma(\overline{\alpha_{i}})] \right)_{1 \leq i \leq g}
$$
whose components are expressed as a computable power series in $s_{1},..., s_{g}$. 
Therefore, the assertion follows from the calculations in Section 2.4. 
\ $\square$ 
\vspace{2ex}

\noindent
{\bf Corollary 3.3.} 
\begin{it} 
Let the notation be as above. 
Then the kernel form  $(g_{{\mathcal R}_{I}})_{f, w}(\tilde{u})$ given in (3.2) is expressed as 
\begin{eqnarray*}
& & 
\left\{ 1 + \sum_{i=1}^{g} \left( 
e^{2 \pi \sqrt{-1} t_{i}} \prod_{k} [y_{k}, x_{k}; \alpha_{i}, \overline{\alpha_{i}}] 
+ e^{-2 \pi \sqrt{-1} t_{i}} \prod_{k} [y_{k}, x_{k}; \alpha_{i}, \overline{\alpha_{i}}] ^{-1} \right) 
\sqrt{s_{i}} \right\} 
\frac{\sqrt{du} \sqrt{d \beta}}{\tilde{u} - \tilde{\beta}} 
\\ 
& & 
+ \ O \left( {\rm max}_{k,l} \{ \sqrt{s_{k}} \sqrt{s_{l}} \} \right),
\end{eqnarray*}
where $\tilde{t} = (t_{1},..., t_{g}) \in {\mathbb R}^{g}$ and 
$\tilde{u} + \tilde{d}(w) = \sum_{k} (y_{k} - x_{k})$.   
\end{it}
\vspace{4ex}

\noindent 
{\bf Data availability} 
\vspace{1ex}

No data was used for the research described in this article. 
\vspace{2ex}

\noindent 
{\bf Acknowledgments} 
\vspace{1ex}

The author would like to thank the reviewers for their helpful comments 
in modifying this manuscript. 
This work is partially supported by the JSPS Grant-in-Aid for 
Scientific Research No. 20K03516.
\vspace{2ex}

\noindent 
{\bf Conflict of interest statement} 
\vspace{1ex}

The author declares that he has no known competing financial interests 
or personal relationships that could have appeared to influence the work 
reported in this article.

\renewcommand{\refname}{\bf References}

\end{document}